\documentclass[conference]{IEEEtran}
\IEEEoverridecommandlockouts
\usepackage{amsmath,amssymb,amsfonts}
\usepackage{algorithmic}
\usepackage{graphicx}
\usepackage{textcomp}
\usepackage{xcolor}
\usepackage{hyperref}

\begin{document}

\title{A Brief Tutorial on Consensus ADMM for Distributed Optimization with Applications in Robotics}

\author{
    \IEEEauthorblockN{Jushan Chen}
    \IEEEauthorblockA{\textit{University of Illinois Urbana-Champaign} \\
    jushanc2@illinois.edu}
}

\maketitle

\begin{abstract}
This paper presents a tutorial on the Consensus Alternating Direction Method of Multipliers (Consensus ADMM) for distributed optimization, with a specific focus on applications in multi-robot systems. In this tutorial, we derive the consensus ADMM algorithm, highlighting its connections to the augmented Lagrangian and primal-dual methods. Finally, we apply Consensus ADMM to an example problem for trajectory optimization of a multi-agent system.
\end{abstract}

\section{Introduction}
In recent years, multi-robot systems have gained significant attention due to their potential in various applications such as surveillance, environmental monitoring, search and rescue, and logistics. Coordination and cooperation among multiple robots, such as drones, are essential to accomplish complex tasks that are beyond the capability of a single agent. One of the key challenges in multi-robot systems is to design algorithms that enable distributed decision-making while ensuring global objectives are met.

Distributed optimization provides a powerful framework for solving large-scale problems in a decentralized manner. It allows multiple agents to optimize their local objectives while cooperating to achieve a global consensus. This is particularly useful in multi-robot systems where communication constraints, scalability, and robustness are critical considerations. For instance, in trajectory planning and collaborative tasks, each robot must compute its own control inputs while considering the actions of its neighbors to avoid collisions and ensure formation maintenance.

The Alternating Direction Method of Multipliers (ADMM) has emerged as an effective optimization algorithm for distributed settings due to its ability to decompose complex problems into simpler subproblems. The Consensus ADMM, a variant of ADMM, is specifically designed for problems where agents need to reach a consensus on shared variables. It combines the advantages of dual decomposition and augmented Lagrangian methods, enabling efficient and scalable solutions.

In this tutorial, we present a derivation of the consensus ADMM algorithm, providing insights into its mathematical foundations based on the augmented Lagrangian and primal-dual methods. We then apply the consensus ADMM algorithm to solve a distributed Model Predictive Control (DMPC) problem in a multi-drone waypoint navigation scenario. This example demonstrates the practical relevance of consensus ADMM in coordinating multiple drones to achieve complicated tasks such as obstacle avoidance.

The rest of the paper is organized as follows. Section .\ref{Sec:background} discusses the mathematical background of distributed optimization. Section .\ref{sec:admm_Derivation} discusses the derivation of the consensus ADMM algorithm. Section .\ref{sec:admm_Application} presents an example problem of applying the consensus ADMM algorithm to a multi-drone waypoint navigation scenario.

\section{Distributed Optimization}
\label{Sec:background}
Distributed optimization deals with solving an optimization problem that is spread across multiple agents. Each agent has access to only part of the global data and cooperates with others to optimize a shared objective function.

The general distributed optimization problem can be formulated as \cite{SOVA_admm}:
\begin{equation}
\begin{array}{ll}
\underset{x_1, \cdots, x_N}{\operatorname{minimize}} & \sum_{i=1}^N g_i\left(\theta_i\right) \\
\text { subject to } & \theta_i=\theta_j \quad \forall j \in \mathcal{N}_i \quad i=1, \cdots, N
\end{array}
\label{eqn:distributed_problem}
\end{equation}
where $\theta_i \in \mathbb{R}^n$ is local state variable of agent $i$, $\theta_j \in \mathbb{R}^n$ is the local state variable of some other agent $j$ that belongs to the neighborhood of agent $i$ (denoted as $\mathcal{N}_i$), and $g_i(\theta_i)$ is the local objective function for agent $i$. A \textit{consensus constraint} between agent $i$ and $j$ is enforced to steer agent $i$'s local state variable $\theta_i$ to agree with the local state variable $\theta_j$ of agent $j$, whenever agent $i$ and $j$ share an edge $\mathcal{E}_{ij}$.

\section{Consensus ADMM Derivation}
\label{sec:admm_Derivation}
The consensus ADMM algorithm is a variant of the classic ADMM algorithm \cite{Boyd2011ADMM}, tailored for distributed optimization problems with consensus constraints. The classic ADMM algorithm prohibits distributed updates of the dual variables, whereas consensus ADMM allows for fully distributed primal-dual updates. In this section, we formulate a derivation of the consensus ADMM algorithm.

\subsection{Augmented Lagrangian}
To solve the constrained optimization problem in Eqn. \ref{eqn:distributed_problem}, we first form the augmented Lagrangian as the following:
\begin{align}
    \mathcal{L}_{\rho}(\theta_i, \theta_{-i}, \lambda_{ij}) &= \sum_{i=1}^{N} g_i(\theta_i) + \sum_{(i,j) \in \mathcal{E}} \Big( \lambda_{ij}^\top (\theta_i - \theta_j) \notag \\
    &\quad + \frac{\rho}{2} \| \theta_i - \theta_j \|^2 \Big)
\end{align}
where $\theta_{-i}$ denotes all agents other than $i$. We denote $\lambda_{ij} \in \mathbb{R}^n$ as the dual variables (Lagrange multipliers) associated with the consensus constraints $\theta_i = \theta_j$ and denote $\rho > 0$ as the penalty parameter for the quadratic penalty term.

Consider the augmented Lagrangian for agent \(i\) in the context of distributed optimization with consensus constraints:
\begin{align}
    \mathcal{L}_\rho(\theta_i, \theta_{-i}, \lambda_i) &= g_i(\theta_i) + \sum_{j \in \mathcal{N}_i} \lambda_{ij}^\top (\theta_i - \theta_j) \notag \\
    &\quad + \frac{\rho}{2} \sum_{j \in \mathcal{N}_i} \|\theta_i - \theta_j\|^2
    \label{eqn:lagrangian_original}
\end{align}
where \( \lambda_{ij} \) are the dual variables associated with the consensus constraints \( \theta_i = \theta_j \). We can simplify this Lagrangian by combining all terms related to agent \(i\). The expression becomes:
\begin{align}
    \mathcal{L}_\rho(\theta_i, \theta_{-i}, \lambda_i) &= g_i(\theta_i) + \sum_{j \in \mathcal{N}_i} \lambda_{ij}^\top \theta_i \notag \\
    &\quad - \sum_{j \in \mathcal{N}_i} \lambda_{ij}^\top \theta_j \notag \\
    &\quad + \frac{\rho}{2} \sum_{j \in \mathcal{N}_i} \|\theta_i - \theta_j\|^2
\end{align}

We now rearrange the Lagrangian to group all the terms involving \(\theta_i\):
\begin{align}
    \mathcal{L}_\rho(\theta_i, \theta_{-i}, \lambda_i) &= g_i(\theta_i) + \left( \sum_{j \in \mathcal{N}_i} \lambda_{ij} \right)^\top \theta_i \notag \\
    &\quad - \sum_{j \in \mathcal{N}_i} \lambda_{ij}^\top \theta_j \notag \\
    &\quad + \frac{\rho}{2} \sum_{j \in \mathcal{N}_i} \|\theta_i - \theta_j\|^2
    \label{eqn:lagrangian_rearranged}
\end{align}

The second term \( \sum_{j \in \mathcal{N}_i} \lambda_{ij}^\top \theta_j \) in Eqn~\eqref{eqn:lagrangian_rearranged} does not depend on \(\theta_i\), so it can be ignored when minimizing with respect to \(\theta_i\). We now expand the quadratic penalty term in Eqn~\eqref{eqn:lagrangian_original}:
\begin{align}
    \frac{\rho}{2} \sum_{j \in \mathcal{N}_i} \|\theta_i - \theta_j\|^2 &= \frac{\rho}{2} \sum_{j \in \mathcal{N}_i} \left( \|\theta_i\|^2 - 2 \theta_i^\top \theta_j + \|\theta_j\|^2 \right)
\end{align}
This term can be grouped with the linear term \( \sum_{j \in \mathcal{N}_i} \lambda_{ij}^\top \theta_i \), resulting in the following simplified expression for the Lagrangian:
\begin{align}
    \mathcal{L}_\rho(\theta_i, \theta_{-i}, \lambda_i) &= g_i(\theta_i) + \left( \sum_{j \in \mathcal{N}_i} \lambda_{ij} - \rho \sum_{j \in \mathcal{N}_i} \theta_j \right)^\top \theta_i \notag \\
    &\quad + \frac{\rho}{2} \|\theta_i\|^2
\end{align}
The primal update for \( \theta_i \) comes from minimizing the Lagrangian with respect to \( \theta_i \), treating \( \lambda_{ij} \) and \( \theta_j \) as fixed. To solve for the primal update, we take the gradient of the Lagrangian with respect to \( \theta_i \) and set it equal to zero:
\[
\nabla_{\theta_i} \mathcal{L}_\rho = \nabla g_i(\theta_i) + \sum_{j \in \mathcal{N}_i} \lambda_{ij} + \rho \sum_{j \in \mathcal{N}_i} (\theta_i - \theta_j) = 0
\]

Solving this for \( \theta_i^{k+1} \) gives:
\begin{align}
\theta_i^{k+1} &= \arg \min_{\theta_i} \Bigg( f_i(\theta_i) 
    + \left( \sum_{j \in \mathcal{N}_i} \lambda_{ij} \right)^\intercal \theta_i \notag \\
    &\quad + \rho \sum_{j \in \mathcal{N}_i} \left\| \theta_i - \frac{\theta_i^k + \theta_j^k}{2} \right\|^2 \Bigg)
\end{align}
\subsection{ADMM Primal-Dual Update Equations}
In the consensus ADMM algorithm, each agent $i$ iteratively updates its primal variables $\theta_i$ by first minimizing the augmented Lagrangian with respect to $\theta_i$ and then updating the dual variables $p_{i}$ using gradient ascent:
\begin{enumerate}
    \item \textbf{Local Primal Update}:
    \begin{equation}
    \begin{aligned}
    \theta_i^{k+1}= & \underset{\theta_i}{\operatorname{argminimize}}\left(g_i\left(\theta_i\right)+\lambda_i^{k \intercal} \theta_i\right. \\
    & \left.+\rho \sum_{j \in \mathcal{N}_i}\left\|\theta_i-\frac{\theta_i^k+\theta_j^k}{2}\right\|_2^2\right)
    \end{aligned}
    \label{eqn:primal_update}
    \end{equation}
    \item \textbf{Dual Variable Update}:
    \begin{equation}
    \lambda_i^{k+1}=\lambda_i^k+\rho \sum_{j \in \mathcal{N}_i}\left(\theta_i^{k+1}-\theta_j^{k+1}\right)
    \label{eqn:dual_update}
    \end{equation}
\end{enumerate}
The primal update~\eqref{eqn:primal_update} minimizes the local objective $g_i(\theta_i)$ while enforcing consensus with the neighboring agents using the quadratic penalty term and the dual variables $\lambda_{ij}$. The dual update~\eqref{eqn:dual_update} adjusts the Lagrange multipliers based on the difference between the current states of neighboring agents. In summary, the consensus ADMM algorithm is an iterative ``negotiation" process where each agent $i$ keeps a \textit{local copy} of the global optimization variable, and it communicates with its neighbor $j$ to solve its local optimization problem.

\section{Application: Multi-Drone Navigation via Distributed MPC}
\label{sec:admm_Application}
We now apply consensus ADMM to a multi-drone waypoint navigation problem formulated as a distributed Model Predictive Control (MPC) problem. Each drone navigates to its goal position by minimizing its objective function while ensuring consensus with neighboring drones to avoid collisions.

\subsection{Problem Description}
We consider a group of $N$ drones connected by a set of edges and vertices $\mathcal{G} =(\mathcal{E},\mathcal{V})$, and each drone $i$ is considered a node on the graph $\mathcal{G}$. The goal of each drone $i$ is to navigate to its final goal position $p_{i}^{ref}\in \mathbb{R}^3$ from its initial position $p_{i}^{0} \in \mathbb{R}^3$ while avoiding collisions with its neighbors. 

\subsection{Drone Dynamics Model}
Each drone is modeled using 12-degree-of-freedom (12-DOF) dynamics, capturing position, velocity, orientation, and angular velocity. We denote the state variable of each drone as:
\begin{align}
    x_i(t) = \begin{bmatrix}
        p_i(t) \\
        v_i(t) \\
        \theta_i(t) \\
        \omega_i(t)
    \end{bmatrix} \in \mathbb{R}^{12}
\end{align}
where $p_i(t) \in \mathbb{R}^3$ is the position, $v_i(t) \in \mathbb{R}^3$ is the velocity, $\theta_i(t) \in \mathbb{R}^3$ is the orientation angles, and $\omega_i(t) \in \mathbb{R}^3$ is the angular velocity. We denote the control input variable of each drone as:
\begin{align}
    u_i(t) = \begin{bmatrix}
        F_i(t) \\
        \tau_i(t)
    \end{bmatrix} \in \mathbb{R}^4
\end{align}
with thrust force denoted as $F_i(t)$ and torques denoted as $\tau_i(t) \in \mathbb{R}^3$. The dynamics of each drone can be denoted as:
\begin{equation}
    \dot{x}_{i}(t) = f_i(x_{i}(t),u_{i}(t))
    \label{eqn:dynamics}
\end{equation}
where $f(\cdot)$ refers to the nonlinear dynamics function of the drone. We assume that each drone has the same dynamics. We consider a finite-horizon optimal control setting for the waypoint navigation scenario, and thus it is necessary to discretize the dynamics function of the drones. We apply the Runge-Kutta 4-th order method to numerical discretize the dynamics equation~\ref{eqn:dynamics}.

\subsection{Distributed MPC Problem Formulation}
\label{sec:mpc_constrained}
We now formulate the distributed MPC problem for the multi-drone waypoint navigation scenario. We denote $N$ as the total number of agents on the graph $\mathcal{G}$. We denote $n$ as the combined dimension of state variables of all agents in $\mathcal{G}$ and similarly denote $m$ as the combined dimension of the input variables of all agents in $\mathcal{G}$. Note that each agent $i$ keeps a \textit{local copy} of the states of all other agents on the graph $\mathcal{G}$. We denote $\tilde{\theta}_{i} \in \mathbb{R}^{(H+1)\times n + H \times m}$ as the local state variable of agent $i$. Intuitively, $\tilde{\theta}_{i}$ refers to the concatenated state and input state variables of all agents in $\mathcal{G}$ with a horizon length of $H$. This definition is adopted for easy adaptation to the consensus ADMM algorithm. We denote the discretized dynamics of the drone as $\bar{f}(\cdot)$. We denote the feasible set of states of agent $i$ as $\mathcal{X}_i$, and similarly denote the feasible set of input of agent $i$ as $\mathcal{U}_i$. Each agent $i$ solves the following optimal control problem:
\begin{align}
    \min_{\{x_{i}(k), u_i(k) \}} \quad & \sum_{k=0}^{H-1} \left( \| x_{i}(k) - x_i^{\text{ref}}(k) \|_Q^2 + \| u_i(k) \|_R^2 \right) \label{eq:local_cost} \\
    \text{subject to} \quad & x_i(k+1) = \bar{f}_i(x_i(k), u_i(k)) \\
    & x_i(0) = x_i^0 \\
    & x_i(k) \in \mathcal{X}_i, \quad u_i(k) \in \mathcal{U}_i \\
    & \| p_i(k) - p_j(k) \| \geq d_{\text{min}}, \quad \forall j \in \mathcal{N}_i
    \label{eqn:local_problem}
\end{align}
where $H$ is the prediction horizon, $Q \succeq 0$, $R \succ 0$ are weight matrices in the quadratic tracking cost, and $d_{\text{min}} > 0$ is the threshold collision avoidance distance.

\subsection{Consensus ADMM Algorithm for Distributed MPC}
We now present the consensus ADMM-MPC problem by reformulating~\eqref{eqn:local_problem} embedding it into the consensus ADMM algorithm. The complete algorithm proceeds as follows:
\begin{enumerate}
    \item \textbf{Local Primal Update}: Each drone solves the following optimization problem:
    \begin{equation}
    \begin{aligned}
    & \tilde{\theta}_i^{t+1}=\underset{\tilde{\theta}_i \in \Theta_i}{\operatorname{argminimize}}\left(g_i\left(\tilde{\theta}_i\right)+p_i^{t \mathcal{\intercal}} \tilde{\theta}_i\right. \\
    &\left.+\rho \sum_{j \in \mathcal{N}_i}\left\|\tilde{\theta}_i-\frac{\tilde{\theta}_i^k+\tilde{\theta}_j^k}{2}\right\|_2^2\right)
    \end{aligned}
    \label{eqn:primal_update}
    \end{equation}
    where $\Theta_i$ is the feasible set of states and inputs for all agents. This step optimizes each drone's local trajectory while considering its neighbors' influence via the consensus terms. The primal update~\eqref{eqn:primal_update} is equivalent to solving~\eqref{eqn:local_problem}, where the quadratic tracking cost and the relevant problem constraints are taken into account.
    
    \item \textbf{Dual Variable Update}: After computing the new trajectories, update the dual variables:
    \begin{equation}
        \tilde{\lambda}_i^{k+1}=\tilde{\lambda}_i^k+\rho \sum_{j \in \mathcal{N}_i}\left(\tilde{\theta}_i^{k+1}-\tilde{\theta}_j^{k+1}\right)
    \end{equation}
    This step enforces consensus between neighboring drones.
    
    \item \textbf{Communication}: Drones exchange their updated positions and states with their neighbors.
    
    \item \textbf{Iteration}: Repeat the above process until convergence such that all the drones have reached their goal positions.
\end{enumerate}

\bibliographystyle{IEEEtran}
\bibliography{citations}

\end{document}